\documentclass[twocolumn]{autart}

\usepackage{amsmath}
\usepackage{amssymb}
\usepackage{hyperref}
\usepackage[english]{babel}

\usepackage{graphicx}
\usepackage{natbib}
\usepackage{lpic}
\usepackage{xcolor}
\usepackage{enumerate}
\usepackage{hyperref}
\usepackage{xcolor}
\usepackage{subfig}
\usepackage{array}
\usepackage{epstopdf}
\usepackage{algorithm}
\usepackage{algpseudocode}
\newcounter{thmcount}

\DeclareMathOperator*{\argmin}{argmin}
\DeclareMathOperator{\co}{Co}
\DeclareMathOperator{\vol}{vol}

\newcommand{\defeq}{:=}

\DeclareMathOperator{\Co}{Co}

\def\RR{\mathbb{R}}
\def\NN{\mathbb{N}}
\def\PP{\mathbb{P}}

\def\X{\mathcal{X}}
\def\U{\mathcal{U}}
\def\S{\mathcal{S}}

\def\W{\mathcal{W}}
\def\M{\mathcal{M}}
\def\probw{p_w}

\def\Id{{\mathcal {{I}}}}

\DeclareMathOperator{\tr}{tr}
\DeclareMathOperator{\vect}{vec}
\DeclareMathOperator*{\minimize}{minimize}
\newcommand{\norm}[1]{\lVert{#1}\rVert}

\newtheorem{theorem}{Theorem}
\newtheorem{lemma}[theorem]{Lemma}
\newtheorem{proposition}[theorem]{Proposition}

\newtheorem{remark}[theorem]{Remark}

\newtheorem{assumption}{Assumption}

\allowdisplaybreaks[4]

\begin{document}
\begin{frontmatter}
\title{Robust Adaptive Model Predictive Control with Persistent Excitation Conditions}

\author[Oxford]{Xiaonan Lu}\ead{wxluxiaonan@hotmail.com}, 
~
\author[Oxford]{Mark Cannon}\ead{mark.cannon@eng.ox.ac.uk}

\address[Oxford]{Department of Engineering Science, University of Oxford, OX1 3PJ, UK}

\begin{abstract}
For constrained linear systems with bounded disturbances and parametric uncertainty, we propose a robust adaptive model predictive control strategy with online parameter estimation.  Constraints enforcing persistently exciting closed loop control actions are introduced for a set-membership parameter identification scheme. The algorithm  requires the online solution of a convex program, satisfies constraints robustly, and ensures recursive feasibility and input-to-state stability. 
Almost sure convergence to the actual system parameters is demonstrated under assumptions on stabilizability, reachability, and tight disturbance bounds.

\end{abstract}
\end{frontmatter}

\thispagestyle{empty}
\pagestyle{plain}

\section{Introduction}
The performance of Model Predictive Control (MPC) relies on an accurate model of the controlled system.
To reduce model uncertainty, adaptive MPC algorithms have been proposed that allow  model parameters to be estimated online. In system identification and adaptive control, persistent excitation (PE) conditions play a key role in establishing convergence of parameter estimates~\citep{Green1986,Shimkin1987}. By incorporating constraints to ensure appropriate PE conditions, a constrained MPC strategy can impose a lower bound on the expected rate of parameter convergence. As a result, adaptive MPC has the potential to estimate system parameters while controlling the system subject to constraints. Various approaches have been proposed~\citep{Mayne2014}, but robust, computationally tractable adaptive MPC remains an open topic under research. 

Adaptive MPC strategies usually have the dual purpose of regulating the system via feedback and providing sufficient excitation for identification of system parameters.  Different adaptive MPC approaches place varying emphasis on these two competing objectives. Some focus on robust constraint satisfaction and stability (e.g.~through constraint tightening~\citep{DiCairano2016}, min-max cost formulations \citep{Adetola2009,Wang2017} or tube MPC \citep{Lorenzen2018,Lu2019}) but omit persistent excitation conditions in the problem formulation.
On the other hand, some approaches consider a nominal MPC problem and force the control law to be persistently exciting, but fail to ensure constraint satisfaction and closed loop system stability \citep{Goodwin1984,Marafioti2014}. 

Other approaches aim to achieve the dual objectives of system regulation and sufficient excitation simultaneously. 
For example, \citet{Weiss2014} select the control input from the Robust Admissible Invariant (RAI) set and balance the conflicting objectives by using an augmented cost function. The resulting control law is more likely to be persistently exciting, but this is not guaranteed. 
\citet{Tanaskovic2014} avoid imposing PE conditions by considering the discrepancy between the nominal and the actual models. However, the proposed algorithm involves solving a nonconvex, infinite-dimensional optimization and can only be simplified for specific examples. 
\citet{Gonzalez2014} separate the two objectives using a dual mode control strategy that injects persistent excitation into the system whenever the state enters a target region for parameter identification.
The proposed algorithm is only applicable to open-loop stable linear systems, and the existence of the target region is example-dependent.
\citet{HernandezVicente2019} propose an algorithm that satisfies PE condition and state and input constraints recursively, but the system model cannot be adapted online. 

In addition, although the importance of persistent excitation conditions have been widely acknowledged in adaptive control literature \citep{Narendra1987}, few strategies incorporate these conditions in a convex optimization formulation. For example, \citet{Marafioti2014} expresses the conditions for persistent excitation as nonconvex quadratic inequalities in terms of the control input. Similarly, \citet{HernandezVicente2019} demonstrate that a PE condition can be satisfied using a periodic solution computed offline, but this solution might not be optimal. Other approaches \citep{Lu2019,Lu2020} use linearization of the PE condition around a reference trajectory to determine sufficient conditions for persistency of excitation, but are unable to ensure closed loop satisfaction of PE conditions through recursively feasible constraints.

In this work we consider linearly constrained linear systems with parametric uncertainty and bounded additive disturbances.
Building on \citet{Lu2020}, we propose an adaptive MPC algorithm that combines set-based parameter identification, robust regulation, and recursively feasible constraints. The algorithm is input-to-state practically stable and provides convergence of parameter estimates almost surely. 
This paper has three main contributions:
\begin{itemize}
\item 
We show that a robust MPC law can be made persistently exciting by adding a random signal to the terminal control law and incorporating additional conditions in the optimization of predicted performance.
\item
We provide a recursively feasible set of conditions for ensuring a persistently exciting control law. These consist of convex constraints in the MPC optimization and a posterior check performed on a nonconvex condition via a set of sampled convex conditions.
\item 
The proposed algorithm  is computationally tractable while achieving the dual objectives of robustly stabilizing the system and providing persistent excitation for online parameter estimation.
\end{itemize}

The remainder of the paper consists of four main sections. Section \ref{sec:PE_linear_feedback} begins by introducing the class of constrained uncertain linear systems under consideration. Then we briefly recap the set-based parameter identification method and introduce the definition of persistent excitation  that is employed in this work. Following this we provide conditions such that the system under a linear feedback law (with or without injected noise) satisfies the required PE condition and hence that the resulting estimated parameter set converges almost surely to the true parameter value. 
Section~\ref{sec:adaptive_MPC} explains the tube-based robust MPC formulation, including the parameterized control law, the initial and terminal conditions, and the cost function. We then define a sequence of PE constraints and show how to convert the resulting nonconvex constraints into convex constraints in the MPC optimization and how to perform a posterior check using sampling. The robust adaptive algorithm is then summarized  and the section concludes with a set of results that demonstrate recursive feasibility, input-to-state practical stability (ISpS) and satisfaction of PE conditions in closed loop operation. Section~\ref{sec:example} provides a numerical example. We compare the proposed algorithm with robust MPC without PE constraints and hence demonstrate that the proposed algorithm results in faster convergence of the estimated parameter set, a higher PE coefficient value and improved performance.

\textit{Notation:} The set of real numbers is denoted $\mathbb{R}$. The non-negative (or positive) integers are denoted $\mathbb{N}_{\geq 0}$ (or $\mathbb{N}_{> 0}$), and $\mathbb{N}_{[p,q]}$ denotes $\{n\in\mathbb{N} : {p \leq n \leq q}\}$.  
The identity matrix is $\Id$. The $i$th element of a vector $a$ is $[a]_i$ and $\norm{a}$ denotes the Euclidean norm.  
The p-dimensional closed unit ball is $\mathbb{B}^p = \{ x\in\mathbb{R}^p : \|x\| \leq 1\}$.
The $i$th row of a matrix $A$ is $[A]_i$, 
$\vect(A)$ is the vector formed by stacking the columns of $A$, 
the matrix inequality $A \geq 0$ applies elementwise, and $A \succeq 0$ (or $A \succ 0$) indicates that $A$ is positive semidefinite (positive definite).  
For $\mathcal{X},\mathcal{Y}\subset\mathbb{R}^n$, $A\mathcal{X} = \{Ax: x\in\mathcal{X}\}$, $\mathcal{X}\oplus\mathcal{Y} = \{x+y: x\in\mathcal{X},\, y\in\mathcal{Y}\}$
and $\Co\{\mathcal{X},\mathcal{Y}\}$ denotes the convex hull of $\mathcal{X}$ and $\mathcal{Y}$.
A polytope is a convex and compact polyhedral subset of $\RR^n$.
The $k$ steps ahead predicted value of a variable $y$ is denoted $y_k$, and the complete notation $y_{k|t}$ indicates the $k$ steps ahead prediction at time $t$.
Probabilities and expectations 
conditioned on the system state $x_t$ are denoted $\mathbb{P}(\cdot |x_t)=\mathbb{P}_t(\cdot)$ and $\mathbb{E}(\cdot|x_t) = \mathbb{E}_t(\cdot)$ respectively, and $\mathbb{P}(\cdot)$, $\mathbb{E}(\cdot)$ are respectively equivalent to 
$\mathbb{P}_0(\cdot)$, $\mathbb{E}_0(\cdot)$.

\section{Persistent excitation conditions} \label{sec:PE_linear_feedback}

In this section we introduce our assumptions on the system model before discussing set-based identification methods. The conditions for persistence of excitation are considered and we show that these are satisfied under linear feedback with mild assumptions. 

\subsection{Problem formulation} \label{subsec:problem_formulation}
The system state $x_t\in \RR^{n_x}$, control input $u_t \in \RR^{n_u}$ and unknown disturbance input $w_t \in\RR^{n_w}$, satisfy
\begin{equation} \label{eq:sys}
x_{t+1} = A(\theta^\ast)x_t +B(\theta^\ast)u_t+ F w_t .
\end{equation}
at all times $t \in \NN_{\geq 0}$. Matrices $A(\theta^\ast)$, $B(\theta^\ast)$ depend on an unknown constant parameter $\theta^\ast\in\mathbb{R}^p$.

\begin{assumption}[Additive disturbance] \label{ass:disturbance}
The disturbance sequence $\{w_t\in\W,\, t \in\NN_{\geq 0}\}$ is independent and identically distributed (i.i.d.), $\mathbb{E} (w_t) = 0$, $\mathbb{E} (w_t w_t^\top) \succeq \epsilon_w \Id$, $\epsilon_w > 0$, and $\W$ is a known convex polyhedral set.
\end{assumption}

\begin{assumption}[Model parameters]\label{ass:param_uncertainty}
(a). $A(\theta), B(\theta)$ are defined in terms of known matrices $A_{i}, B_{i}$, ${i \in\NN_{[0,p]}}$:
\[
(A(\theta),B(\theta)) = (A_0,B_0) +\smash{\sum_{i = 1}^{p} (A_i,B_i)} [\theta]_i 
, \ \ \forall \theta\in\Theta_0 .
\]
(b). $\Theta_0$ is a known polytopic set containing $\theta^\ast$:
\[
\theta^\ast\in\Theta_0 = \{ \theta : M_\Theta \theta \leq \mu_{0} \} 
= \co \{ \theta_0^{(1)}, \dots, \theta_0^{(m)} \} .
\]
(c). The pair $\bigl(A(\theta^\ast),B(\theta^\ast)\bigr)$ is reachable.
\\
(d). $(A(\theta),B(\theta))=(A(\theta^\ast),B(\theta^\ast))$ if and only if $\theta=\theta^\ast$.
\end{assumption}

\subsection{Set-based parameter identification}\label{sec:parameter_set_estimation}

Parameter estimation methods include recursive least squares~\citep{Heirung2017}, comparison sets~\citep{Aswani2013}, set-membership identification~\citep{Tanaskovic2014,Lorenzen2018} and neural network training~\citep{Akpan2011,Reese2016}. Here we use a set-membership approach to enable robust satisfaction of constraints. Set-based parameter identification was proposed in~\citep{Chisci1998,Veres1999}  and it was shown in~\cite{Lu2020} that the estimated parameter set converges to the true parameter value with probability~1 if the associated regressor is persistently exciting (PE). 

At times $t\in\mathbb{N}_{>0}$ we use observations of the state $x_t$ to determine a set $\Delta_t$ of consistent model parameters, known as the unfalsified parameter set. This is combined with $\Theta_{t-1}$ to construct a new parameter set estimate $\Theta_t$.

The model (\ref{eq:sys}) can be rewritten as
\[
x_{t+1} = \Phi(x_t, u_t) \theta^\ast + \phi(x_t, u_t)+ F w_t
\]
where $\Phi_t $ and $\phi_t $ are known at time $t$ and are defined by 
\begin{align}
\Phi_t & =\Phi(x_t,u_t)  = \begin{bmatrix} A_1 x_t+B_1 u_t &  \cdots & A_p x_t + B_p u_t\end{bmatrix} 
\label{eq:Phik}
\\
\phi_t & =\phi(x_t,u_t) = A_0 x_{t}+B_0 u_{t} . 
\end{align}
Given $x_t$, $x_{t-1}$, $u_{t-1}$ and the disturbance set $\W$, the unfalsified parameter set at time $t$ is given by
\begin{equation}\label{eq:unfalsified_set_def}
\Delta_t = \{ \theta : x_t -A(\theta)x_{t-1} -B(\theta)u_{t-1} \in  F\W\}
\end{equation}
The parameter set $\Theta_t$ may be updated using $\Delta_t$ by various methods, including minimal~\citep{Chisci1998}, fixed-complexity~\citep{Lorenzen2018}, and limited-complexity~\citep{Tanaskovic2014} update laws. In each case, $\Theta_t$ is non-increasing and $\Theta_t\subseteq \Theta_{t-1}$ for all~$t\in\mathbb{N}_{> 0}$.

For a fixed-complexity parameter set update law, 
the parameter set estimate $\Theta_t$ is defined as
$
\Theta_t  = \Theta(\mu_t) = \{ \theta : M_{\Theta} \theta \leq \mu_{t} \}
$
where  $M_\Theta\in\RR^{r\times p}$ is an \textit{a priori} chosen matrix and {$\mu_t \in \RR^r$} is determined 
so that $\Theta_t$ is the smallest set containing the intersection of $\Theta_{t-1}$ and the unfalsified sets
$\Delta_j$ for $j=t-N_\mu+1,\ldots,t$, 
\begin{equation} \label{eq:fixed_set_update}
\mu_{t} := \min_{\mu\in\RR^r} {\vol} \bigl(\Theta(\mu)\bigr)
\ \ \text{s.t.} \!
\smash{\bigcap_{j=t-N_\mu+1}^t}  \!\! \Delta_j \cap \Theta_{t-1} \! \subseteq \Theta(\mu)
\end{equation}
where $\Delta_j := \RR^{p}$ for $j\leq0$ and $N_\mu$ is the parameter update window length. Note that $\mu_t$ can be computed by solving a set of linear programs.

We briefly recap the definition of persistent excitation (PE). 
The regressor $\Phi_t$ in (\ref{eq:Phik}) is persistently exciting if there exist a horizon $N_u$ and a scalar $\epsilon_\Phi>0$ such that
\begin{equation}\label{eq:pe_general}
\sum_{k=t}^{t+N_u-1} \Phi_k^\top \Phi_k  \succeq \epsilon_\Phi \Id
\end{equation}
for all $t \in \NN_{\geq 0}$. 
In the current work however, we define persistent excitation using the expectation:
\begin{equation}\label{eq:pe_expectation}
\sum_{k=t}^{t+N_u-1} \mathbb{E} \bigl\{ \Phi_k^\top \Phi_k\bigr\} \succeq \epsilon_\Phi \Id ,
\end{equation}
which is required to hold for some $\epsilon_\Phi > 0$, at an infinite number of discrete time instants $t \in\mathbb{N}_{\geq0}$. 
We refer to the interval $\mathbb{N}_{[ t, t+N_u-1]}$ as a PE window.

\begin{assumption}[Tight disturbance bound]\label{ass:probw}
For any point $w^0$ on the boundary of $\W$
and any $\epsilon > 0$, the disturbance sequence $\{w_0,w_1,\ldots\}$ satisfies
$\PP \bigl\{ \norm{ w_t - w^0 } < \epsilon \bigr\} \geq \probw (\epsilon) > 0$,
for all $t\in\NN_{\geq 0}$.
\end{assumption}

If the regressor in~(\ref{eq:Phik}) satisfies the PE condition~(\ref{eq:pe_general}), 
then under Assumptions~\ref{ass:disturbance} and \ref{ass:probw}
the parameter set estimate $\Theta_t$ under a fixed-complexity update law with $N_\mu \geq N_u$
converges to $\{ \theta^\ast\}$ with probability~1.
This is shown in~\citet[][Theorem~3 and Corollary~3]{Lu2020}.
Here we use a version of this result that applies to~(\ref{eq:pe_expectation}).

\begin{lemma}\label{lem:convergence}
Under Assumptions~\ref{ass:disturbance} and \ref{ass:probw},
if $\Phi_t$ satisfies the PE condition~(\ref{eq:pe_expectation}),
then the minimal complexity parameter set estimate $\Theta_t$ 
with $N_\mu\geq N_u$
converges to $\{\theta^\ast\}$ as $t\to\infty$ with probability~1.
\end{lemma}

\textbf{Proof.}\hspace{1ex}If~(\ref{eq:pe_expectation}) holds,
then there exists 
(with probability~1) an infinite sequence $\{\kappa_i, \, i\in\NN_{\geq 0}\}$ satisfying, for some $\epsilon_\Phi > 0$ and any given $\theta\in\mathbb{R}^p$,
\[
\sum_{k=\kappa_i}^{\kappa_i+N_u-1} \| \Phi_k\theta\|^2   \geq \epsilon_\Phi \lVert\theta\rVert^2.
\]
Choose $\kappa_i$ so that $\kappa_{i+1}\geq \kappa_i+N_u-1$ and therefore $\Theta_{\kappa_i+N_u-1} \supseteq \Theta_{\kappa_{i+1}}$. Then we can prove the Lemma using Theorem~3 and Corollary~3 of~\citet{Lu2020} with $t$ replaced by $\kappa_t+N_u-1$, since these results imply that any $\theta\in\Theta_0$ such that $\theta\neq \theta^\ast$ is necessarily excluded from $\Theta_{\kappa_i}$ with a probability that converges to $1$ as $i\to\infty$.
It follows that 
$\Theta_t\to\{\theta^\ast\}$ as $t\to\infty$ with probability~1.
\qed

Assumption~\ref{ass:disturbance} on the disturbance sequence $\{w_t,t\in\mathbb{N}_{\geq 0}\}$ is common in robust MPC and set-membership identification formulations. 
Clearly, if the disturbance sequence contains temporal correlation that is representable as an i.i.d.\ disturbance sequence driving a known linear filter, then this can be incorporated in the dynamics (\ref{eq:sys}) without violating Assumption~\ref{ass:disturbance}.
In addition, the assumption of zero-mean disturbances can be relaxed without affecting the convergence result (Lemma~\ref{lem:convergence}). For example, by redefining 
the model parameters $A(\theta),B(\theta),F$ and $\theta$, the mean value of the disturbance can be estimated with the other unknown model parameters via the parameter set update law.
Assumption~\ref{ass:probw} may be difficult to verify, but we note that this assumption can be relaxed at the expense of some residual uncertainty in the parameter set estimate~\citep[][Sec.~5.3]{Lu2020}. 
We also note that, although $x_t$ is assumed to be exactly known in the definition (\ref{eq:unfalsified_set_def}) of the unfalsified set $\Delta_t$, this assumption can be relaxed to allow the use of noisy measurements or state estimates \citep[see e.g.][Sec.~5.4]{Lu2020}.

In the next section we  show that the PE condition~(\ref{eq:pe_expectation}) can be satisfied under a given linear feedback law, which provides useful insight into the properties of predicted control laws. 

\subsection{Regressor under linear feedback law}

Consider the system (\ref{eq:sys}) under the action of a linear feedback law $u_t=Kx_t$. To simplify notation we define $A_K(\theta)=A(\theta) + B(\theta)K$ and $A_{K,i} = A_i + B_i K$, $i\in\mathbb{N}_{[0,p]}$.

\begin{assumption}[Stability]\label{ass:quad_stable}
For $z_t\in\RR^{n_x}$ and $t\in\NN_{\geq 0}$, $z_{t+1} \in\Co\{ A_K(\theta) z_t,$ $\theta\in\Theta_0\}$ is quadratically stable.
\end{assumption}

\begin{assumption}\label{ass:reachability}
 The gain $K$ is chosen so that:\\ 
 (a). $A_K(\theta) =  A_K(\theta^\ast)$ if and only if $\theta=\theta^\ast$.
 \\
 (b). The pair $\bigl( A_K(\theta^\ast), F\bigr)$ is reachable.
 \end{assumption}

Assumption~\ref{ass:quad_stable} is equivalent to feasibility of the linear matrix inequalities (LMIs): $P - \smash{A_K(\theta_0^{(j)})^\top P A_K(\theta_0^{(j)})} \succ 0$, $j\in\mathbb{N}_{[1,m]}$, $P = P^\top\succ 0$ \citep[e.g.][Theorem 6.21]{blanchini2008}.

Assumption~\ref{ass:reachability}(b) is more difficult to verify because $\theta^\ast$ is unknown. However, a sufficient condition is that $( A_K(\theta_t), F)$ is reachable for all sequences $\{\theta_0,\theta_1,\ldots\}$  such that $\theta_t \in \Theta_0$ for all $t$, which can be verified by checking the feasibility of a set of LMIs defined by the vertices of $\Theta_0$.

\begin{theorem} \label{thm:PE_vector}
If $N_u> n_x$, then under Assumptions \ref{ass:disturbance}, \ref{ass:param_uncertainty} and \ref{ass:reachability},
the regressor $\Phi_t$ in (\ref{eq:Phik}) with $u_t = Kx_t$ satisfies~(\ref{eq:pe_expectation}) for some $\epsilon_\Phi > 0$ and all $t\in\NN_{\geq 0}$,
\end{theorem}

\textbf{Proof.}\hspace{1ex}Let $N_u> n_x$. Under the feedback law $u_t = K x_t$, the system~(\ref{eq:sys}) becomes $x_{t+1} = A_K(\theta^\ast) x_t + F w_t$, and the definition of the regressor $\Phi_t$ yields, for all $\theta\in \RR^p$,
\begin{equation}\label{eq:expansion_phi}
\begin{aligned}
& \sum_{k=t}^{t+N_u-1} \theta^\top  \mathbb{E}_t \bigl[ \Phi(x_k,Kx_k)^\top \Phi(x_k,Kx_k) \bigr] \theta \\
& \quad = 
\tr \Bigl\{
\sum_{i=1}^p A_{K,i}^\top \, [\theta]_i
\sum_{j=1}^p A_{K,j} \, [\theta]_j 
\sum_{k=t}^{t+N_u-1} \mathbb{E}_t(x_k x_k^\top)
\Bigr\} .
\end{aligned}
\end{equation}
But $w_k$ is i.i.d. and $\mathbb{E} (w_k w_k^\top)  \succeq \epsilon_w \Id$ (Assumption~\ref{ass:disturbance}), so
$\mathbb{E}_t ( x_k x_k^\top) 
\succeq \epsilon_w\sum_{j=0}^{k-t-1} \bigl[ A_K(\theta^\ast)\bigr]^j F F^\top \bigl[A_K(\theta^\ast)^\top\bigr]^j $
for all $k> t$.
In addition, the pair $\bigl(A_K(\theta^\ast),F\bigr)$ is reachable (Assumption~\ref{ass:reachability}), so a positive scalar $\epsilon_{\mathcal{R}}$ exists such that $\sum_{j=0}^{n_x-1} \bigl[ A_K(\theta^\ast)\bigr]^j F F^\top \bigl[ A_K(\theta^\ast)^\top\bigr]^j \succeq \epsilon_{\mathcal{R}} \Id$. It follows that
$\mathbb{E}_t ( x_k x_k^\top)\succeq \epsilon_w \epsilon_{\mathcal{R}} \Id$
for all $k \geq t+n_x$ and hence
\[
 \sum_{k=t}^{t+N_u-1} \mathbb{E}_t(x_k x_k^\top) 
\succeq \epsilon_w \epsilon_{\mathcal{F}}\Id
 \]
where $\epsilon_{\mathcal{F}} = (N_u - n_x) \epsilon_{\mathcal{R}} $. 
Therefore (\ref{eq:expansion_phi}) becomes
\begin{align*}
&\sum_{k=t}^{t+N_u-1} \theta^\top \mathbb{E}_t\bigl[ \Phi(x_k,Kx_k)^\top \Phi(x_k,Kx_k) \bigr] \theta
\\
& \quad \geq  \epsilon_w \epsilon_{\mathcal{F}} 
\bigl\|\vect \Bigl( \smash{\sum_{i=1}^p A_{K,i}} \, [\theta]_i \Bigr) \bigr\|^2
\geq \sigma^2 \epsilon_w \epsilon_{\mathcal{F}} \norm{\theta}^2
\end{align*}
where we have used the identities $\tr(A^\top A) =  \norm{\vect(A)}^2$ and $\vect \bigl( \sum_{i=1}^p A_{K,i} \, [\theta]_i \bigr)  = \bigl[  \vect(A_{K,1}) \ \ \cdots \ \ \vect( A_{K,p} ) \bigr] \theta $. Here $\sigma$ is 
the minimum singular value of the matrix $[ \vect(A_{K,1}) \ \cdots \ \vect( A_{K,p} ) ]$, which 
is non-zero due to Assumption~\ref{ass:reachability}. It follows that~(\ref{eq:pe_expectation}) holds for all $t\geq 0$ with $\epsilon_{\Phi} = \sigma^2 \epsilon_w \epsilon_{\mathcal{F}} >0$ for any $N_u> n_x$.
\qed

Lemma~\ref{lem:convergence} and Theorem~\ref{thm:PE_vector} imply that the minimal and fixed-complexity parameter set estimates under linear feedback converge to $\{\theta^\ast\}$ with probability~1 if Assumptions \ref{ass:disturbance}-\ref{ass:reachability} hold.

\subsection{Regressor under linear feedback with injected noise}\label{sec:with_noise}
Theorem~\ref{thm:PE_vector} links the convergence of parameter set estimates to properties of the disturbance sequence (Assumptions~\ref{ass:disturbance} and \ref{ass:probw}) and $(A_i, B_i)$ and $K$ (Assumptions~\ref{ass:param_uncertainty}  and \ref{ass:reachability}). 
 To relax these conditions, specifically Assumptions~\ref{ass:reachability}(a) and (b), consider a feedback law with injected noise: 
\begin{equation}\label{eq:control2}
u_t = Kx_t+ s_t.
\end{equation}
Here $s_t \in \RR^{n_u}$ is a random variable with a known probability distribution that provides additional excitation to encourage convergence of parameter estimates. 

\begin{assumption}\label{ass:noise_s}
The sequence $\{s_t\in\S, \, t\in \NN_{\geq 0}\}$ is i.i.d.\ with $\mathbb{E}(s_t)= 0$, $\mathbb{E}(s_t s_t^\top)\succeq \epsilon_s \Id$, $\epsilon_s>0$, 
 $s_t$ is independent of $x_t$ and $w_t$, and $\S$ is a known polytopic set. 
\end{assumption}

Similarly to Theorem \ref{thm:PE_vector}, we next give a bound on the expectation of the PE condition (\ref{eq:pe_general})
under feedback law~(\ref{eq:control2}).

\begin{theorem}\label{thm:PE_vector2}
If $N_u > n_x$, then under Assumptions \ref{ass:disturbance}, \ref{ass:param_uncertainty} and \ref{ass:noise_s},
the regressor $\Phi_t$ in (\ref{eq:Phik}) with $u_t = Kx_t +s_t$ satisfies (\ref{eq:pe_expectation}) for some $\epsilon_{\Phi} >0$ and all $t\in\NN_{\geq 0}$.
\end{theorem}

\textbf{Proof.}\hspace{1ex}Let $N_u > n_x$ be given. The control law (\ref{eq:control2}) implies $x_{t+1} = A_K(\theta^\ast) x_t + B(\theta^\ast) s_t + Fw_t$, and since $\bigl(A_K(\theta^\ast),B(\theta^\ast)\bigr)$ is reachable (Assumption \ref{ass:param_uncertainty}(c)) we have
\[
\mathbb{E}_t(x_k x_k^{\top}) \succeq \epsilon_s \! \sum_{j=0}^{k-t-1} A_K(\theta^\ast)^{j\!} B(\theta^\ast) B(\theta^\ast)^{\!\top\!\!}
{A_k(\theta^\ast)^{j\!}}^\top\!\! \succeq \epsilon_s \epsilon_{\mathcal{R}}
\mathcal{I} 
\]
for all $k\geq t+n_x$, where $\epsilon_{\mathcal{R}} > 0$.

Hence Assumption~\ref{ass:noise_s} implies
\[
\sum_{k=t}^{t+N_u-1}\!\!\!\mathbb{E}_t\biggl\{ \begin{bmatrix} x_k \\ s_k\end{bmatrix} \! \begin{bmatrix} x_k^{\top\!} & s_k^{\top\!}\end{bmatrix} \biggr\}
\succeq \!\sum_{k=0}^{N_u-1} \!\! \begin{bmatrix} \mathbb{E}_t[x_k x_k^\top]\!\! & 0 \\ 0 & \!\!\epsilon_s \Id \end{bmatrix} 
\succeq \epsilon_s\epsilon_{\mathcal{F}} \Id
\]
where $\epsilon_{\mathcal{F}} = \min\{(N_u - n_x) \epsilon_{\mathcal{R}}, 1\}$. 
Furthermore, since
$A_ix_k+B_iu_k = A_{K,i}x_k + B_is_k = [ A_{K,i} \ B_i  ][ x_k^\top \ s_k^\top]^\top$,
an argument similar to the proof of Theorem~\ref{thm:PE_vector} yields
\begin{align*}
& \sum_{k=t}^{t+N_u-1} \theta^\top \mathbb{E}_t \bigl[\Phi(x_k,Kx_k+s_k)^\top \Phi(x_k,Kx_k+s_k) \bigr] \theta 
\\
&\quad \geq \epsilon_s\epsilon_{\mathcal{F}} \tr\Bigl\{ 
\sum_{i = 1}^p\begin{bmatrix} A_{K,i}^\top \\ B_i^\top \end{bmatrix} [\theta]_i 
\sum_{j = 1}^p\begin{bmatrix} A_{K,j} & B_j \end{bmatrix} [\theta]_j 
\Bigr\}
\\
 &\quad \geq \sigma^2 \epsilon_s\epsilon_{\mathcal{F}} \norm{\theta}^2
\end{align*}
for all $\theta\in\RR^p$, where $\sigma$ is the minimum singular value of the matrix whose $i$th column is $[ \vect(A_{K,i})^\top \ \vect(B_i)^\top ]^\top$ for $i\in\mathbb{N}_{[1,p]}$ and $\sigma>0$ by Assumption~\ref{ass:param_uncertainty}(d). Hence (\ref{eq:pe_expectation}) holds with {$\epsilon_\Phi = \sigma^2\epsilon_s\epsilon_{\mathcal{F}} >0$}.
\qed

\begin{remark}
In comparison to Theorem \ref{thm:PE_vector}, Theorem \ref{thm:PE_vector2} is less restrictive because it does not rely on Assumption~\ref{ass:reachability} and 
Assumption~\ref{ass:noise_s} can be easily satisfied because $s_t$ is user-defined. 
\end{remark}

It follows from Lemma~\ref{lem:convergence}  and Theorem~\ref{thm:PE_vector2} that Assumptions \ref{ass:disturbance}-\ref{ass:probw} and \ref{ass:noise_s}  ensure that  the minimal and fixed-complexity parameter set estimates under $u_t = Kx_t+s_t$ converge to $\{\theta^\ast\}$ with probability~1.

\section{Adaptive robust model predictive control}\label{sec:adaptive_MPC}

The noise input $s_t$ added to the feedback law in (\ref{eq:control2}) can cause poor tracking performance and may violate constraints on the states and control inputs of the system~(\ref{eq:sys}). However, a receding horizon control law incorporating injected noise can avoid these undesirable effects while exploiting the PE properties it provides.
Consider a predicted control law parameterized at time $t$ in terms of decision variables $\mathbf{v}_t = \{v_{0|t},\ldots,v_{N-1|t}\}$:
\begin{align} \label{eq:input_law}
u_{k|t} =\begin{cases}
Kx_{k|t}+v_{k|t}, & k\in \mathbb{N}_{[0, N-1]}\\
Kx_{k|t} +s_{k|t}, & k \in \mathbb{N}_{[N,N+N_u-1]}
\end{cases}
\end{align}
where $N$ is the conventional MPC prediction horizon \citep[see e.g.][]{Kouvaritakis2015}. The sequence $\mathbf{s}_t=\{\smash{s_{k|t}},\, k\in\mathbb{N}_{[0,N+N_u-1]}\}$ is a random sequence such that $\mathbf{s}_{t}$ is independent of $\mathbf{s}_{q}$ for all $t\neq q$, $t,q\in\mathbb{N}_{\geq 0}$ and $\mathbf{s}_{t}$ satisfies Assumption~\ref{ass:noise_s} (i.i.d.\ with $\mathbb{E}(s_{k|t}) = 0$, $\mathbb{E}(s_{k|t}\smash{s_{k|t}^\top}) \succeq \epsilon_s\Id \succ 0$, $s_{k|t} \in \mathcal{S}$, $\mathbf{s}_{t}$ independent of $x_{t}$ and $w_{t}$).
  The feedback gain $K$ is assumed to be stabilizing (Assumption~\ref{ass:quad_stable}).

We assume linear state and control input constraints
\begin{equation} \label{eq:input_state_constraint}
x_{k|t} \in \mathcal{X}, \quad u_{k|t} \in \U, \quad 
 \forall k \in \NN_{[0,N+N_u-1]} ,
\end{equation}
where $\X$, $\U$ are given polytopes.
To enforce these constraints we define a terminal set $\X_T$ satisfying 
 \begin{align}
 &\X_T \subseteq\X, \quad K\X_T \oplus \mathcal{S} \subseteq \U ,
 \label{eq:terminala}
 \\
 &A_K(\theta) \X \oplus B(\theta)\mathcal{S} \oplus F \W \subseteq \X_T \quad
 \forall \theta \in \Theta_t .
 \label{eq:terminalb}
 \end{align}

\subsection{Tube MPC formulation}

To ensure satisfaction of constraints (\ref{eq:input_state_constraint}),
we construct a sequence of sets denoted 
$\mathbf{X}_t = \{\X_{k|t},\, k\in\mathbb{N}_{[0,N+N_u-1]}\}$,
 satisfying, for all $\theta\in\Theta_t$ and $k\in\mathbb{N}_{[0,N-1]}$,
 \begin{align}
& \X_{k|t} \subseteq \X, \quad K \X_{k|t} \oplus \{ v_{k|t}\} \subseteq \U ,
\label{eq:tube_xu_constraint}
\\
& A_K(\theta)\X_{k|t} \oplus \{ B(\theta)v_{k|t}\} \oplus F\W \subseteq \X_{k+1|t} ,
\label{eq:recurrence_constraint}
\end{align}
with initial and terminal conditions, 
\begin{align}
\X_{0|t} &= \{x_{t} \} ,
\label{eq:initial_condition}
\\
\X_{k|t} & = \X_T, \ k \in \mathbb{N}_{[N,N+N_u-1]} ,
\label{eq:terminal_constraint}
\end{align}
where $\X_T$ satisfies~(\ref{eq:terminala})-(\ref{eq:terminalb}).

The tube $\mathbf{X}_t$ is an online optimization variable, and, depending on the computational resources available, different tube formulations may be employed. Here we assume homothetic~\citep{Lorenzen2018} or fixed-complexity polytopic tubes~\citep{Lu2019} that have convenient vertex representations:
 $\X_{k|t} = \smash{\Co\{x_{k|t}^{(1)},\ldots, x_{k|t}^{(\nu)}\}}$,
where $\nu$ is a fixed number of vertices. 

Various objective functions have been proposed for adaptive MPC. Examples include min-max piecewise linear~\cite{Lu2019}, min-max quadratic~\cite{Adetola2009}, and nominal quadratic~\cite{Lorenzen2018} cost functions.
Here we consider a nominal cost defined for a given nominal parameter vector $\bar{\theta}_t\in\Theta_t$ by
\[
J(x_t,{\bf v}_t,\bar{\theta}_t) = \sum_{k = 0}^{N-1} l(\bar{x}_{k|t},Kx_{k|t}+v_{k|t}) + V_{N|t}(\bar{x}_{N|t})
\]
with $\bar{x}_{0|t}=x_t$, $\bar{x}_{k+1|t} = A_K(\bar{\theta}_t) \bar{x}_{k|t} + B(\bar{\theta}_t) v_{k|t}$, $k\in\mathbb{N}_{[0,N-1]}$. 
The stage cost $l(\cdot,\cdot)$ and terminal cost $V_{N|t}(\cdot)$ are assumed to be positive definite quadratic functions
satisfying, for given $\bar{\theta}_t\in\Theta_0$ and all $x\in\mathbb{R}^{n_x}$,
 \begin{equation}\label{eq:terminal_cost}
V_{N|t}(x) = V_{N|t}\bigl(A_K(\bar{\theta}_t)x\bigr) +l(x, Kx).
 \end{equation}
In the remainder of the paper we use $(\mathbf{v}_t^o, \mathbf{X}_t^o)$ to denote the optimal decision variables $(\mathbf{v}_t,\mathbf{X}_t)$ for the problem of minimizing $J(x_t,\mathbf{v}_t,\bar{\theta}_t)$ subject to constraints (\ref{eq:tube_xu_constraint})-(\ref{eq:terminal_constraint}) and additional constraints to ensure the PE condition~(\ref{eq:pe_expectation}), as we discuss next in Section~\ref{subsec:PE_condition_formulation}.

\subsection{PE condition} \label{subsec:PE_condition_formulation}

To define a recursively feasible set of constraints, we construct a series of $N_{pe}$ overlapping PE windows extending from the past and across the prediction horizon, where
\[N_{pe} = \begin{cases} N+t ,  & t < N_u-1 \\ 
N+N_u-1 , &  t \geq N_u-1 .
\end{cases}
\]Consider the PE conditions defined at time $t\in\mathbb{N}_{> 0}$ for given $\mathbf{X}_t$, $\mathbf{v}_t$, $\mathbf{s}_t$,  
 and all $\kappa \in \NN_{[N-N_{pe}, N-1]}$ by
 \begin{equation}\label{eq:pe_condition_window}
 \sum_{k=\kappa}^{\kappa+N_u-1} \!\!\! \Phi (x_{k|t}, Kx_{k|t}+q_{k|t})^{\!\top} \Phi(x_{k|t}, Kx_{k|t}+q_{k|t})  \succeq \beta_{\kappa|t} \Id
 \end{equation}
for all $x_{k|t}\in\X_{k|t}$, $k\in\mathbb{N}_{[\kappa,\kappa+N_u-1]}$, where $q_{k|t} := v_{k|t}$  if $k < N$ and
 $q_{k|t} := s_{k|t}$ if $k\geq N$. The maximum value of $\beta_{\kappa|t}$ satisfying (\ref{eq:pe_condition_window}) defines a PE coefficient corresponding to a PE window starting at time $t+\kappa$. With a slight abuse of notation, we allow $\kappa < 0$ and thus include past states and control inputs in this definition by defining $\mathcal{X}_{k|t} := \{x_{k+t}\}$, {$Kx_{k|t}+v_{k|t} := u_{k+t}$} whenever $k < 0$.
 
The PE conditions in (\ref{eq:pe_condition_window}) are nonconvex in  $(\mathbf{v}_{t}, \mathbf{X}_t)$, and hence are unsuitable as  constraints in an online MPC optimization. However, following \cite{Lu2020}, 
 we can linearize these conditions around a reference trajectory $(\hat{\mathbf{x}}_t, \hat{\mathbf{u}}_t) = \{(\hat{x}_{k|t}, \hat{u}_{k|t}), \, k\in\mathbb{N}_{[0,N+N_u-1]}\}$.
For a given nominal parameter vector $\bar{\theta}_t$ and sequences $\hat{\mathbf v}_t$, $\mathbf{s}_t$, the reference trajectory is defined by
\begin{subequations} \label{eq:reference_xvu}
\begin{align}
& \hat{x}_{0|t} = x_t \label{eq:reference_x0}\\
 & \hat{x}_{k+1|t} = A(\bar{\theta}_{t})\hat{x}_{k|t} +B(\bar{\theta}_t)\hat{u}_{k|t} , \,k\in\mathbb{N}_{[0,N+N_u-2]} \!\!
 \\
 & \hat{u}_{k|t} = \begin{cases}
 K\hat{x}_{k|t}+\hat{v}_{k|t} , & k\in\mathbb{N}_{[0,N-1]} \\
 K\hat{x}_{k|t}+ {s}_{k|t} , & k\in\mathbb{N}_{[N,N+N_u-1]} .
 \end{cases}
 \label{eq:reference_uk}
\end{align}
\end{subequations}
We define the sequence $\hat{\mathbf v}_t$ for $t\in\mathbb{N}_{> 0}$ using the solution of the MPC optimization at time $t-1$, denoted $\mathbf{v}^o_{t-1}
 =\{ v_{0|t-1}^o, \ldots, v_{N-1|t-1}^o\}$, 
 and $s_{N-1|t}$:
 \begin{equation}\label{eq:reference_vk}
 \hat{v}_{k|t} = \begin{cases}
 v^o_{k+1|t-1} ,  & k\in\mathbb{N}_{[0,N-2]} \\
 s_{N-1|t}  , & k=N-1.
 \end{cases}
 \end{equation}
A nominal parameter vector $\bar{\theta}_t\in\Theta_t$ can be obtained for example by projecting an RLS parameter estimate onto $\Theta_t$ \citep{Lorenzen2018} or by projecting $\bar{\theta}_{t-1}$ onto $\Theta_t$,
 \begin{equation} \label{eq:nominal_theta_update}
 \bar{\theta}_{t} = \argmin_{\theta\in\Theta_t} \|\bar{\theta}_{t-1} - \theta\| .
 \end{equation}

Linearizing (\ref{eq:pe_condition_window}) by neglecting quadratic terms in the decision variables $(\mathbf{v}_t,\mathbf{X}_t)$ yields a set of LMIs in $\mathbf{v}_t$, the vertices, ${x_{k|t}^{(j)}}$, of $\mathbf{X}_t$, and additional optimization variables $\boldsymbol{\beta}_t' = \{\beta_{\kappa | t}', \, \kappa \in \NN_{[N-N_{pe}, N-1]}\}$,
 given for $\kappa \in \NN_{[N-N_{pe}, N-1]}$ by
 \begin{equation}\label{eq:PE_windows}
 \sum_{k=\kappa}^{-1} \Phi_{k+t}^\top \Phi_{k+t} 
 + \sum_{k = \max\{0,\kappa\}}^{\kappa+N_u-1} M_{k|t} 
 \succeq  {\beta}_{\kappa|t}' \Id 
 \end{equation}
 where $M_{k|t} = M_{k|t}^\top \in\RR^{p\times p}$, $k\in\NN_{[0,N+N_u-2]}$ satisfies 
 \begin{align*}
 M_{k|t} &\preceq
 \hat{\Phi}_{k|t}^{\top} \hat{\Phi}_{k|t} 
 +
 \hat{\Phi}_{k|t}^\top \Phi(x^{(j)}_{k|t} - \hat{x}_{k|t}, Kx^{(j)}_{k|t} +q_{k|t}- \hat{u}_{k|t}) \\
 &\quad   +\Phi(x^{(j)}_{k|t} - \hat{x}_{k|t}, Kx^{(j)}_{k|t} +q_{k|t} - \hat{u}_{k|t}) ^\top  \hat{\Phi}_{k|t} ,
 \end{align*}
 for all $j\in\mathbb{N}_{[1,\nu]}$, with $\hat{\Phi}_{k|t} = \Phi(\hat{x}_{k|t},\hat{u}_{k|t})$.
 
To increase the probability that the MPC law satisfies the PE condition~(\ref{eq:pe_expectation}), we 
include (\ref{eq:PE_windows}) in the MPC optimization at times $t>0$ with the following constraints on $\boldsymbol{\beta}_t'$
  \begin{equation}\label{eq:beta_bound_moving}
\smash{ \beta_{\kappa|t}' \geq \hat{\beta}'_{\kappa|t} , \ \ \forall \kappa \in \NN_{[N-N_{pe}, N-1]},}
 \end{equation}
where $\smash{\hat{\beta}_{\kappa|t}^\prime}$ is a lower bound on the maximum value of $\smash{\beta_{\kappa|t}^\prime}$ satisfying (\ref{eq:PE_windows}) for $(\mathbf{v}_t,\mathbf{X}_t)$ satisfying (\ref{eq:tube_xu_constraint})-(\ref{eq:terminal_constraint}). Thus we determine
 $\smash{\hat{\boldsymbol{\beta}}\mbox{}'_{t}} = \{\hat{\beta}_{\kappa | t}', \, \kappa \in \NN_{[N-N_{pe}, N-1]}\}$
 by finding
   $\smash{\hat{\beta}_{\kappa|t}^\prime}$ in (\ref{eq:beta_bound_moving}) as the solution of
 \begin{equation}\label{eq:reference_beta_moving}
 \hat{\beta}_{\kappa|t}^\prime := \max_{\beta_{\kappa|t}^\prime\in\mathbb{R}} \beta_{\kappa|t}^\prime\ \ \text{s.t.}\ \ (\ref{eq:PE_windows}) 
 \end{equation}
with $\smash{q_{k|t} := \hat{v}_{k|t}}$ if $k < N$ and $\smash{q_{k|t} := s_{k|t}}$ if $k\geq N$, and with $x_{k|t}^{(j)} := x_{k+1|t-1}^{(j)\,o}$, $\forall j\in\mathbb{N}_{[1,\nu]}$, $k >0$, where $\smash{x_{k|t-1}^{(j)\,o}}$ denotes a vertex of $\mathbf{X}_{t-1}^o$.

{To ensure that (\ref{eq:pe_expectation}) holds with non-zero probability for all $t$, we propose a check on the solution of the online MPC optimization based on the nonconvex condition}
\begin{equation}\label{eq:beta_recursive}
\beta_{\kappa|t}\geq \hat{\beta}_{\kappa|t},\qquad \forall \kappa \in \NN_{[N-N_{pe}, N-1]},
\end{equation}
where $\beta_{\kappa|t}$ is defined in (\ref{eq:pe_condition_window}) and $\hat{\beta}_{\kappa|t}$ has the same definition except that $q_{k|t} \defeq \hat{v}_{k|t}$ if $k<N$.

Condition (\ref{eq:beta_recursive}) is difficult to check directly because (\ref{eq:pe_condition_window}) is nonconvex in $x_{k|t}$ and $q_{k|t}$. Therefore we check instead an approximate condition using sampling. 
For each $k \in \NN_{[0,N+N_u-1]}$, we draw $N_s$ samples $\smash[b]{\{{x}_{k|t}^{(1)},\ldots,{x}_{k|t}^{(N_s)}\}}$ uniformly distributed on $\X_{k|t}$ and compute the corresponding regressor products defined by
\[
\Psi^{(i)}_{k|t}(\mathbf{v}_t) = \Phi(x^{(i)}_{k|t},Kx^{(i)}_{k|t}+v_{k|t})^\top \Phi(x^{(i)}_{k|t},Kx^{(i)}_{k|t}+v_{k|t})
\]
for $i \in \NN_{[1,N_s]}$.
We then impose the condition $ \smash{\beta^s_{\kappa|t}\geq \hat{\beta}^s_{\kappa|t}}$ with sampled states, where $ \beta^s_{\kappa|t}$ and $\hat{\beta}^s_{\kappa|t}$ are defined as
 \begin{align*}
\beta^s_{\kappa|t}&:= \max_\beta \beta \ \text{s.t.}  \min_{\mathbf{i} \in \NN_{[1,N_s]} \times \cdots \times \NN_{[1,N_s]}} \!\!\!\!
\sum_{k=\kappa}^{\kappa+N_u-1}\!\!\! \Psi^{(i_k)}_{k|t}(\mathbf{v}^o_t)  \succeq \beta \Id 
\\
\hat{\beta}^s_{\kappa|t}&:= \max_\beta \beta \ \text{s.t.} 
\min_{\mathbf{i} \in \NN_{[1,N_s]} \times \cdots \times \NN_{[1,N_s]}} \!\!\!\!
\sum_{k=\kappa}^{\kappa+N_u-1}\!\!\! \Psi^{(i_k)}_{k|t}(\hat{\mathbf{v}}_t)  \succeq \beta \Id 
 \end{align*}
where $\mathbf{i} = (i_\kappa , \ldots, i_{\kappa+N_u-1})$.
Since the samples $x^{(i)}_{k|t}$ are uniformly distributed over $\X_{k|t}$, and the product $\smash{\Psi^{(i)}_{k|t}(\mathbf{v}_t)} $ is continuous in $\smash{x^{(i)}_{k|t}}$, we can infer that (\ref{eq:beta_recursive}) is satisfied with a probability of at least $p_s$, where $p_s>0$, whenever $ \beta^s_{\kappa|t}\geq \hat{\beta}^s_{\kappa|t}$.
The adaptive MPC algorithm is summarized in Algorithm~\ref{alg:noise_s}.


\begin{algorithm}
\caption{Adaptive MPC with PE constraints}\label{alg:noise_s}
\textbf{\\At} $t=0$:
Choose $\Theta_0$ and $\mathcal{S}$. Determine $K$ satisfying Assumption~\ref{ass:quad_stable} and $\X_T$, $V_{N|0}$ satisfying (\ref{eq:terminala})-(\ref{eq:terminalb}) and~(\ref{eq:terminal_cost}). Define $N$, $N_u$, $N_w$, $N_\theta$. Obtain $\bar{\theta}_0$ and $x_0$
and compute the solution, $(\mathbf{v}^o_0, \mathbf{X}^o_0)$, of the quadratic program (QP):
\begin{equation} \label{eq:online_mpc_offline}
\begin{aligned}
& \mathcal{P}_{0}: \ \minimize_{\mathbf{v}_0,\mathbf{X}_0} \ J(x_0,\mathbf{v}_0,\bar{\theta}_0)
\ \ \text{s.t.\ (\ref{eq:tube_xu_constraint})-(\ref{eq:terminal_constraint})}.
\end{aligned}
\end{equation}
Apply the control input $u_0 = K x_0 + v_{0|0}^o$.

\textbf{At times} $t=1,2,\ldots$:
\begin{enumerate}[(a).]
\item Obtain the current state $x_t$.
\item Update $\bar{\theta}_t$, $\Theta_{t} $ and $V_{N|t}$ using(\ref{eq:fixed_set_update}), (\ref{eq:nominal_theta_update}) and (\ref{eq:terminal_cost}).
\item 
Generate the noise sequence $\mathbf{s}_t$ and compute $\hat{\mathbf{x}}_t $, $\hat{\mathbf{u}}_t $, $\hat{\mathbf{v}}_t$, $\smash{\hat{\boldsymbol{\beta}}\mbox{}'_{t}}$ using (\ref{eq:reference_xvu}), (\ref{eq:reference_vk}) and (\ref{eq:reference_beta_moving}). 
\item 
Compute the solution $(\mathbf{v}^o_t, \mathbf{X}_t^o)$ of the semidefinite program (SDP):
\begin{equation}\label{eq:online_mpc}
\begin{aligned}
\mathcal{P}_{>0}:\   &\minimize_{\mathbf{v_t},\mathbf{X}_t,\boldsymbol{\beta}_t} \ 
J(x_t,\mathbf{v}_t,\bar{\theta}_t)
\\
&\ \text{s.t.\ (\ref{eq:terminala})-(\ref{eq:terminal_constraint}), (\ref{eq:PE_windows}) and (\ref{eq:beta_bound_moving})}.
\end{aligned}
\end{equation}
\item 
Generate $N_s$ samples of \smash{$\mathcal{X}_{k|t}$} and $N_s$ samples of \smash{$\mathcal{X}^o_{k+1|t-1}$} for each $k\in\smash{\mathbb{N}_{[\kappa,\kappa+N_u-1]}}$ and compute {$\beta^s_{\kappa|t}$} and {$\hat{\beta}^s_{\kappa|t}$} for $\kappa\in\NN_{[N-N_{pe},N-1]}$. If \smash{$\beta^s_{\kappa|t} < \hat{\beta}^s_{\kappa|t}$} for any $\kappa$, set $\mathbf{v}_t^o := \hat{\mathbf{v}}_t$ and $\mathbf{X}_t^o := \bigl\{\{x_t\},\X^o_{2|t-1}, \ldots,$ $\smash[t]{\X^o_{N-1|t-1},\X_T, \ldots, \X_T \bigr\}}$.
\item
Apply the control input $u_t = K x_t + v_{0|t}^o$
\end{enumerate}
\end{algorithm}

\begin{remark}
PE constraints are applied to $N_{pe}$ windows at each time step, and each constraint (\ref{eq:PE_windows})  is a LMI. Therefore the computational load of online optimization initially increases over time and levels off at time $t = N+N_u-1$.
\end{remark}

\begin{remark}
Alternative parameter set and nominal parameter update laws can be used in Algorithm~\ref{alg:noise_s} instead of (\ref{eq:fixed_set_update}) and (\ref{eq:nominal_theta_update}), provided that $\Theta_{t+1} \subseteq \Theta_t$ and $\bar{\theta}_t \in \Theta_t$.
\end{remark}

\subsection{Recursive Feasibility and Input-to-State Stability} \label{sec:properties1}
\begin{proposition}[Recursive feasibility]\label{prop:feas_moving}
The online MPC optimization in Algorithm \ref{alg:noise_s}
is feasible at all times $t \geq 1$ if (\ref{eq:online_mpc}) is feasible at $t = 1$ and $\Theta_{t}\subseteq\Theta_{t-1}$ for all $t$.
\end{proposition}

\begin{proposition}[Stability]\label{prop:stability}
Let Assumptions \ref{ass:disturbance}, \ref{ass:param_uncertainty} and \ref{ass:quad_stable} hold. Then the system (\ref{eq:sys}) under Algorithm~\ref{alg:noise_s} is input-to-state practically stable (ISpS)~\citep[][Def.~6]{Limon2008} in the set of initial conditions $x_0$ for which $\mathcal{P}_0$ is feasible.
\end{proposition}
\textbf{Proof.}\hspace{1ex}For $\mathbf{s}_t$, $\mathbf{v}^o_{t}$, $\hat{\mathbf{v}}_t$ in Algorithm~\ref{alg:noise_s} steps (c) and (d) let $\zeta^o_{k|t} = v_{k|t}^o - s_{k|t}$ and $\smash{\hat{\zeta}}_{k|t} = \smash{\hat{v}}_{k|t} - s_{k|t}$ for all $k\in\mathbb{N}_{[0,N-1]}$.

Since $\mathcal{P}_{>0}$ is necessarily feasible if $\mathcal{P}_0$ is feasible, input-to-state practical stability follows from {Theorem~1 and Corollary~1} of \cite{Lu2020}. This can be shown using the proofs of \citep[Theorem~1 and Corollary~1]{Lu2020} with $v_{k|t}$ and $\hat{v}_{k|t}$ replaced  respectively by $\zeta_{k|t}$ and $\smash{\hat{\zeta}_{k|t}}$, and with disturbance input $w_{t}$ replaced by $Fw_t + B(\theta^\ast) s_{0|t}$.
 \qed
 
Proposition~\ref{prop:stability} implies~\citep{Limon2008,Lu2020} the existence of a $\mathcal{K}\mathcal{L}$-function $\eta$ and $\mathcal{K}$-functions $\psi$, $\xi$ satisfying, for all $x_0$ such that $\mathcal{P}_0$ is feasible, the bound
 \begin{align*}
 \| x_t\| &\leq \eta(\|x_0\|, t) + \psi\bigl(\max_{\tau\in\NN_{[0,t-1]}} \|Fw_\tau + B(\theta^\ast)s_{0|\tau}\| \bigr)
 \\
 &\quad + \xi\bigl( \max_{\tau\in\NN_{[0,t-1]}} \| \bar{\theta}_\tau - \theta^\ast \|\bigr).
 \end{align*}

\subsection{Closed loop satisfaction of the PE condition} \label{sec:properties2}
We define the closed loop PE coefficient $\epsilon_t$ at any time $t\geq N_u-1$ as the maximum value of $\epsilon_t$ satisfying
\begin{equation} \label{eq:epsilon_t}
 \sum_{k=t-N_u+1}^{t} \Phi (x_k, u_k)^{\top} \Phi(x_{k}, u_{k})  \succeq \epsilon_t \Id .
\end{equation}
To ensure closed loop satisfaction of the PE condition, a more restrictive version of Assumption~\ref{ass:param_uncertainty}(c) is needed.
\begin{assumption}\label{ass:sys_matrices2}
The pair $\bigl(A_K(\theta_t), B(\theta_t)\bigr)$ is reachable for all sequences $\{\theta_0,\theta_1, \ldots\}$ such that $\theta_t \in \Theta_0$ for all $t$. 
\end{assumption}
This reachability assumption can be verified by checking 
whether a set of LMIs defined by the vertices of $\Theta_0$ is feasible.

\begin{theorem} \label{thm:pe_closed_loop}
Under Assumptions \ref{ass:disturbance}, \ref{ass:param_uncertainty}, \ref{ass:noise_s} and \ref{ass:sys_matrices2},
the closed loop PE coefficient $\epsilon_t$ in (\ref{eq:epsilon_t}) under Algorithm~\ref{alg:noise_s} satisfies
\begin{equation}
\mathbb{E}( \epsilon_t ) \succeq \epsilon_\Phi
\end{equation}
for some $\epsilon_\Phi >0$ with non-zero probability, for all $t \geq N+N_u-1$.
\end{theorem}

\textbf{Proof.}\hspace{1ex}Given the optimal control 
$\mathbf{v}_t^o$ at time $t$,  the predicted PE coefficient $\beta_{\kappa|t}$ in (\ref{eq:pe_condition_window}) can be rewritten as
\begin{align}
&\beta_{\kappa|t} =  \min_{x_{k|t}\in \X_{k|t}}  \max_{\beta} \beta \label{eq:beta_n} \\
&\text{s.t. }
\sum_{k = \kappa}^{\kappa+N_u -1} \Phi(x_{k},Kx_{k}+q_{k|t})^\top \Phi(x_{k},Kx_{k}+q_{k|t}) \succeq \beta \Id\nonumber
\end{align}
where $q_{k|t} := v^o_{k|t}$ if $k<N$ and $q_{k|t}:=s_{k|t}$ if $k\geq N$.

Using the reference input $\hat{\bf v}_{t}$ defined in (\ref{eq:reference_vk}), the reference PE coefficient $\hat{\beta}_{\kappa-1|t+1}$ in (\ref{eq:reference_beta_moving}) of the same PE window is
\begin{align}
& \hat{\beta}_{\kappa-1|t+1} = \min_{z_k \in \mathcal{Z}_k} \max_{\beta} \beta\label{eq:beta_z} \\
&\text{s.t.} \sum_{k = \kappa}^{\kappa+N_u -1} \Phi(z_{k},Kz_{k}+{q}_{k|t})^\top \Phi(z_{k},Kz_{k}+{q}_{k|t}) \succeq \beta \Id \nonumber 
\end{align}
where $\mathcal{Z}_{k} = \X_{k-1|t+1}$ and 
$q_{k|t} := \hat{v}^o_{k|t}$ if $k<N-1$ and $q_{k|t}:=s_{k-1|t+1}$ if $k\geq N-1$.

But $\Theta_{t+1} \subseteq \Theta_t$ and the constraints (\ref{eq:initial_condition}) imply 
$\mathcal{Z}_k \subseteq \X_{k|t}$, for all $k \in \NN_{[1, N+N_u-1]}$. 
Given that the injected noise sequences $\mathbf{s}_{t+1}$ and $\mathbf{s}_{t}$ are i.i.d with zero mean and the regressor $\Phi(x,u)$ is affine in $(x,u)$, we can deduce that 
$\mathbb{E}_{t}\bigl[\hat{\beta}_{\kappa-1|t+1} \bigr] \geq \mathbb{E}_t\bigl[  \beta_{\kappa|t}\bigr]$.
Combining this with the bounds enforced by step (e) in Algorithm \ref{alg:noise_s}, we have
\begin{equation}\label{eq:exp_beta_chain}
\mathbb{E}_{t}\bigl[\beta_{\kappa-1|t+1} \bigr] 
\geq
\mathbb{E}_{t}\bigl[\hat{\beta}_{\kappa-1|t+1} \bigr] 
\geq
\mathbb{E}_t\bigl[  \beta_{\kappa|t}\bigr]
\end{equation}
with probability $p_s>0$.

Next we derive a lower bound on $\mathbb{E}_t[ \beta_{N|t}]$. 
Given Assumption~\ref{ass:sys_matrices2}, there exists a scalar $\epsilon_\mathcal{B}>0$ satisfying $\sum_{j= 0}^{n_x-1} [A_K(\theta)]^j B(\theta) B(\theta)^\top [{A_K(\theta)^\top}]^j \succeq \epsilon_\mathcal{B}\Id$, and hence 
$ \sum_{k = N}^{N+N_u -1} \mathbb{E}_T [x_{k|t} x_{k|t}^\top]  \succeq (N_u - n_x)\epsilon_{\mathcal{B}} \mathbb{E}(s_k s_k^\top)
$.
The predicted PE coefficient $\beta_{N|t}$ therefore satisfies
\begin{align*}
  \mathbb{E}_t&[ \beta_{N|t}] 
\overset{(\text{E1})}{=}  \min_{{\theta}\in \mathbb{B}^p} \min_{x_{k|t}\in \X_{k|t}} \Bigg\{
\sum_{k = N}^{N+N_u -1} \mathbb{E}_t \Bigg[ 
s_k^\top B({\theta})^\top B({\theta})s_k
    \\
    & \hspace{30mm}+\Big( A_K({\theta})x_{k|t}\Big)^\top \Big( A_K({\theta})x_{k|t}\Big) 
    \Bigg] \Bigg\}
\\
 \overset{(\text{I2})}{ \geq} &\min_{{\theta}\in \mathbb{B}^p}
  \text{tr}\Bigg\{  \Big( A_K({\theta})\Big)^\top \Big(A_K({\theta})\Big)(N_u - n_x)\epsilon_{\mathcal{B}} \mathbb{E}\Big[s_k s_k^\top\Big]
 \\
 &\hspace{30mm}+B({\theta})^\top B({\theta})\sum_{k = N}^{N+N_u -1}\mathbb{E} \Big[s_k s_k^\top\Big]
 \Bigg\} 
 \\
 \overset{(\text{I3})}{\geq}&  \min_{{\theta}\in \mathbb{B}^p} \Big\{ \bigl\|\text{vec}\Big( \sum_{i=1}^p \begin{bmatrix}\epsilon'_{\mathcal{B}} A_{K,i} &\sqrt{N_u} B_i \end{bmatrix}  [{\theta}]_i\Big)\bigr\|^2 \epsilon_s 
 \Big\}
\\
\overset{(\text{I4})}{\geq} &  \epsilon_s \sigma^2
\end{align*}
where $\epsilon'_{\mathcal{B}} = \sqrt{(N_u - n_x)\epsilon_{\mathcal{B}}}$ and $\sigma$  is the minimum singular value of the matrix whose $i$th column is 
$[\epsilon'_{\mathcal{B}}\vect(A_{K,i})^\top \ \sqrt{N_u}\vect(B_i)^\top ]^\top$ for $i\in\mathbb{N}_{[1,p]}$.
Here  equality (E1) holds because $x_{k|t}$ and the injected noise $s_{k|t}$ are independent and $\mathbb{E}[s_{k|t}] = 0$, inequality (I2) holds by substitution,
inequality (I3) follows from Assumptions \ref{ass:param_uncertainty} and \ref{ass:noise_s}, and $\sigma> 0$ by Assumption \ref{ass:param_uncertainty}.

Applying (\ref{eq:exp_beta_chain}) 
for $\kappa = N,N-1,\ldots,N-N_{pe}+1$,
and noting that 
$\epsilon_{t} = {\beta}_{N-N_{pe}|t}$, we therefore have
\[
\mathbb{E}_{t}\bigl[ \epsilon_{t+{N_{pe}}}   \bigr]
\geq   \mathbb{E}_t\bigl[ \beta_{N|t} \bigr] \geq \epsilon_\Phi
\]
with probability $\smash{p_s^{N_{pe}}}>0$, where $\epsilon_{\Phi} =  \epsilon_s\sigma^2   > 0$.
\qed

If Assumptions \ref{ass:disturbance}-\ref{ass:quad_stable} and \ref{ass:noise_s} hold, then Lemma~\ref{lem:convergence} and Theorem~\ref{thm:pe_closed_loop} imply that the parameter set estimate $\Theta_t$ converges under Algorithm~\ref{alg:noise_s} to $\{\theta^\ast\}$ with probability~1.

\section{Numerical example} \label{sec:example}

In this section, we use the example in~\cite{Lorenzen2018}
to illustrate that the proposed adaptive MPC algorithms can produce persistently exciting input signals while satisfying the system input and state constraints, and providing performance close to the optimal performance for MPC without PE constraints. 

Consider the second-order uncertain linear system~\citep{Lorenzen2018} defined by~(\ref{eq:sys}) with
\begin{gather*}
A_0 = \begin{bmatrix} 0.5 & 0.2 \\ -0.1 & 0.6 \end{bmatrix}, \ 
A_1 = \begin{bmatrix} 0.042 & 0 \\ 0.072 & 0.03 \end{bmatrix}, \\
A_2 = \begin{bmatrix} 0.015 & 0.019 \\ 0.009 & 0.035 \end{bmatrix},  \
A_3 = \begin{bmatrix} 0 & 0 \\ 0 & 0 \end{bmatrix}, \\
B_0 = \begin{bmatrix} 0 \\ 0.5 \end{bmatrix}, \ 
B_1 = \begin{bmatrix} 0 \\ 0 \end{bmatrix}, \ 
B_2  = \begin{bmatrix} 0 \\ 0 \end{bmatrix}, \ 
B_3 = \begin{bmatrix} 0.040 \\ 0.054 \end{bmatrix} ,
\end{gather*}
$\theta^* = [ 0.8 \ \  0.2 \ \  {-0.5}]^{\top\!}$ and $K=[0.017 \ \ {-0.41}]$. The initial parameter set estimate is $\Theta_0 = \{ \theta : \| \theta \|_\infty \leq 1 \}$,  and for all $t\geq 0$, $\Theta_t$ is a {hyperrectangle}, with $\M_\theta = [\Id \ \ {-\Id}]^\top$.
The bounded disturbance sequence $\{w_0,w_1,\ldots\}$ follows a truncated normal distribution with zero mean, covariance $\sigma_w = 0.06$, and bounds $w_t \in \W = \{w: |w |\leq 0.05\}$. The injected noise sequence $\{s_0,s_1,\ldots\}$ follows a uniform  distribution with zero mean and $s_t \in \S = \{s: |s |\leq 0.005\}$. 
The state and input constraints are $[x_t]_2 \geq -0.3$ and $u_t\leq 1$.  
The MPC prediction horizon and PE window length are $N = 10$ and $N_u = 3$ respectively. 
For these parameters we obtain $\epsilon_s = 8.33\times10^{-6}$, $\sigma = 1.27\times10^{-3}$ and $\epsilon_\Phi = 1.35\times10^{-11}$ in Theorem~\ref{thm:pe_closed_loop}.

For comparison with Algorithm~\ref{alg:noise_s} we define a robust adaptive MPC law (Algorithm~\ref{alg:QP}) with online parameter set update but without PE constraints in the MPC optimization.
\begin{algorithm}
\caption{Robust Adaptive MPC} \label{alg:QP}
Modify Algorithm \ref{alg:noise_s} by omitting the constraints 
 (\ref{eq:PE_windows}) and (\ref{eq:beta_bound_moving}) from problem $\mathcal{P}_{>0}$ and omitting step (e).
\end{algorithm}

All simulations were performed in Matlab on a 2.9 GHz Intel Core i7-10700F processor, and the online MPC optimization was formulated using Yalmip~\citep{Yalmip}. Mosek~\citep{Mosek}  and Gurobi~\citep{Gurobi} were used to solve SDP and QP problems respectively. For the purposes of comparison, identical disturbance sequences and the same parameter update laws were employed in each simulation.

Table~\ref{table:algorithmCs11} compares Algorithms~\ref{alg:noise_s},~\ref{alg:QP}, and the linear feedback law with injected noise ($u_t=Kx_t+s_t$), in terms of the parameter set volume (denoted $\lvert \Theta_t \rvert$), empirical mean closed loop PE coefficient ($\hat{\epsilon}$), and average solver run-time.
Here the observed average closed loop PE coefficient is defined as 
\[
 \hat{\epsilon} = \frac{1}{N_{av}N_w}\sum_{i = 1}^{N_{w}} \sum_{k=0}^{N_{av}-1} \epsilon_k
\]
with $N_{av} = 500$ time steps and $N_w = 30$ disturbance sequences. 
Algorithm~\ref{alg:noise_s} outperforms Algorithm~\ref{alg:QP} in terms of the convergence rate of its estimated parameter set, which is consistent with the higher value of $\hat{\epsilon}$ observed for Alg.~\ref{alg:noise_s}. The linear feedback law provides a value of $\hat{\epsilon}$ similar to Alg.~\ref{alg:noise_s}, however this control law does not guarantee constraint satisfaction because the initial state is outside the terminal region.

\renewcommand{\arraystretch}{1.2}
\begin{table}[h]
\caption{Parameter set volume, PE coefficient, and computation}
\label{table:algorithmCs11}
\centerline{\begin{tabular}{m{10em}| m{3.5em} m{3.5em}  m{3.5em}}
\hline
& $Kx_t+s_t$& Alg.~\ref{alg:noise_s} &  Alg.~\ref{alg:QP}\\
\hline
$\lvert\Theta_{100}\rvert/\lvert\Theta_{0}\rvert$ $(\%)$ & 22.82  & 20.83& 22.18  \\
$\lvert\Theta_{200}\rvert/\lvert\Theta_{0}\rvert$ $(\%)$ & 21.47  & 19.93& 21.00  \\
$\lvert\Theta_{500}\rvert/\lvert\Theta_{0}\rvert$ $(\%)$ & 17.24  & 15.54& 17.54  \\
 $\hat{\epsilon}$ ($\times 10^{-7}$) &  $1.042$ & $1.044$ & $0.547$ \\
Step (c) solver time (s) & -  & 0.080  & - \\
Step (d) solver time (s) & -  & 0.027  & 0.075 \\
Step (e) solver time (s) & -  & 3.82  & - \\
 \hline
 \end{tabular}}
\end{table}

The average computation time for the QP in step (d) of Alg.~\ref{alg:QP} exceeds that of the SDP in step (d) of Alg.~\ref{alg:noise_s}. However, Alg.~\ref{alg:noise_s} requires an additional $0.08$\,s per time step to compute $\smash{\hat{\boldsymbol{\beta}}\mbox{}_t^\prime}$ in step (c) and $3.82$\,s to perform step (e) with $N_s = 30$. 
Note that $\smash{\beta_{\kappa|t} < \hat{\beta}_{\kappa|t}^s}$  occurred in step (e) of Alg.~\ref{alg:noise_s} at only $20\%$ of time-steps. This is consistent with the observation that, 
due to the presence of the constraints (\ref{eq:PE_windows}) and (\ref{eq:beta_bound_moving}) in problem $\mathcal{P}_{>0}$, 
Alg.~\ref{alg:noise_s} provides a higher value of $\hat{\epsilon}$ and faster parameter convergence than Alg.~\ref{alg:QP} even if the time-consuming step (e) is omitted from Alg.~\ref{alg:noise_s}, 
although the guarantee of Theorem~\ref{thm:pe_closed_loop} then no longer applies.

\begin{figure}[h]
    \centering
    \includegraphics[scale=0.6, trim = 0mm 0mm 0mm 7.5mm, clip=true]{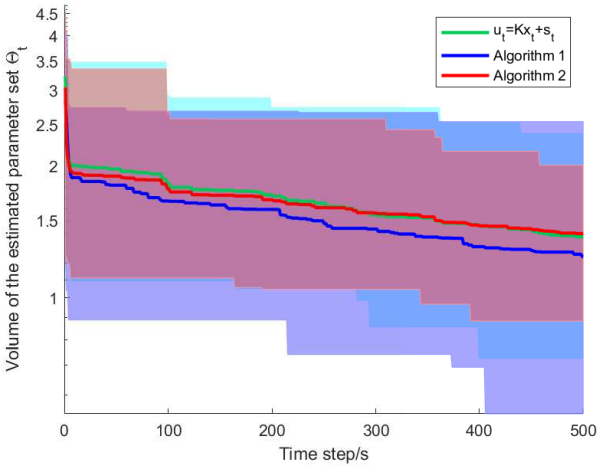}
    \vspace{-6mm}
    \caption{Parameter convergence: Alg.~\ref{alg:noise_s},~\ref{alg:QP} and $u_t=Kx_t+s_t$}
    \label{fig:volume_vs_QP2}
\end{figure}

\begin{figure}[h]
    \centering
    \includegraphics[scale=0.6, trim = 0mm 0mm 0mm 8mm, clip=true]{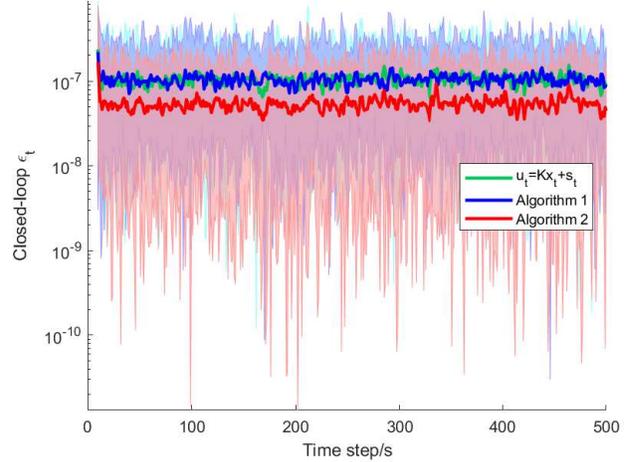}
    \vspace{-6mm}
    \caption{PE coefficient $\epsilon_t$: Alg.~\ref{alg:noise_s}, \ref{alg:QP} and $u_t=Kx_t+s_t$}
    \label{fig:beta_vs_QP2}
\end{figure}

Figure \ref{fig:volume_vs_QP2} compares the volumes of estimated parameter sets for Algorithms~\ref{alg:noise_s},~\ref{alg:QP} and linear feedback with injected noise. For Alg.~\ref{alg:noise_s} (blue) with PE constraints, the volume decreases faster than for Alg.~\ref{alg:QP} (red) without PE constraints or linear feedback with injected noise (green).
Figure \ref{fig:beta_vs_QP2} shows the closed loop value of  $\epsilon_t$ obtained using Algorithms \ref{alg:noise_s},~\ref{alg:QP} and linear feedback with injected noise. As expected from the empirical mean $\hat{\epsilon}$ (Table~1), Alg.~\ref{alg:noise_s} is observed to inject more excitation into the system than Alg.~\ref{alg:QP}, and this explains its faster convergence rate.
The smaller parameter set volume brings 2.2\% improvement in closed loop cost for Alg.~\ref{alg:noise_s} over Alg.~\ref{alg:QP} with identical initial conditions and disturbance sequences.
 
\section{Conclusions}
The proposed robust adaptive MPC algorithm exploits the persistence of excitation properties of linear feedback laws with injected noise while ensuring satisfaction of constraints on states and control inputs, and optimal tracking/regulation performance.
  This is achieved by incorporating recursively feasible constraints into the online MPC optimization which, when combined with a posterior sample-based check, ensures persistence of excitation. 
  Guarantees of input-to-state practical stability and asymptotic parameter set convergence are provided.
 
\bibliographystyle{agsm}
\footnotesize{\bibliography{sample1}}

\end{document}